\documentclass[reqno,12pt]{amsart}
\usepackage{amssymb}
\usepackage{amsfonts}
\usepackage{amsmath}
\usepackage{graphicx}
\usepackage{amscd,amsthm}

\newtheorem{theorem}{Theorem}
\theoremstyle{plain}
\newtheorem{acknowledgement}{Acknowledgement}

\numberwithin{equation}{section}

\input{tcilatex}

\begin{document}
\author{}
\title{}
\maketitle

\begin{center}
\bigskip

\thispagestyle{empty} \pagestyle{myheadings} 
\markboth{\bf Burak Kurt}{\bf
On the degenerate poly-Frobenius-Genocchi polynomials of complex variables
}

\textbf{\Large On The Degenerate Poly-Frobenius-Genocchi Polynomials Of
Complex Variables}

\bigskip

\textbf{Burak Kurt} \bigskip

\medskip

Akdeniz University, Mathematics of Department, Antalya TR-07058,

\qquad Turkey \\[0pt]

E-mail\textbf{: burakkurt@akdeniz.edu.tr}

\medskip

\bigskip

\bigskip

\textbf{{\large {Abstract}}}\medskip
\end{center}

The main aim of this paper is to define and investigate a new class of the
degenerate poly-Frobenius-Genocchi polynomials with the help of the
polyexponential functions. In this paper, we define the degenerate
poly-Frobenius-Genocchi polynomials of complex variables arising from the
modified polyexponential functions. We derive explicit expressions for these
polynomials and numbers. Also, we obtain implicit relations involving these
polynomials and some other special numbers and polynomials.

\begin{quotation}
\ \ 
\end{quotation}



\section{Introduction and Notation}

Throughout this paper, $%
\mathbb{N}
$ denotes the set of \ natural numbers, $%
\mathbb{N}
_{0}$ denotes the set of \ nonnegative integers, $%
\mathbb{R}
$ denotes the set of real numbers and $%
\mathbb{C}
$ denotes the set of complex numbers. We begin by introducing the following
definitions and notations (\cite{1}-\cite{18}).

The Frobenius-Euler polynomials $H_{n}\left( x;u\right) $ are defined by (%
\cite{1}-\cite{18});%
\begin{equation}
\sum\limits_{n=0}^{\infty }H_{n}\left( x;u\right) \frac{t^{n}}{n!}=\frac{1-u%
}{e^{t}-u}e^{xt}\text{,}  \label{1}
\end{equation}%
where $u\neq 1$ and $e^{t}\neq u$.

When $x=0$, $H_{n}\left( u\right) :=H_{n}\left( 0;u\right) $ are called the
Frobenius-Euler numbers.

The Genocchi polynomials are defined by (\cite{11}, \cite{12}, \cite{14})%
\begin{equation}
\sum\limits_{n=0}^{\infty }G_{n}\left( x\right) \frac{t^{n}}{n!}=\frac{2t}{%
e^{t}+1}e^{xt}\text{, }\left\vert t\right\vert <\pi \text{.}  \label{2}
\end{equation}%
When $x=0$, $G_{n}\left( 0\right) :=G_{n}$ are called the Genocchi numbers.

The Frobenius-Genocchi polynomials are defined by \cite{18}%
\begin{equation}
\sum\limits_{n=0}^{\infty }FG_{n}\left( x,u\right) \frac{t^{n}}{n!}=\frac{%
\left( 1-u\right) t}{e^{t}-u}e^{xt}\text{.}  \label{3}
\end{equation}

For $u=-1$, $FG_{n}\left( x,-1\right) =G_{n}\left( x\right) $ and $x=0$, $%
FG_{n}\left( u\right) :=FG_{n}\left( 0,u\right) $ are called the
Frobenius-Genocchi numbers.

The degenerate exponential function is defined by (\cite{3}-\cite{11}) with $%
\lambda \in 
\mathbb{R}
\backslash \left\{ 0\right\} $%
\begin{equation}
e_{\lambda }^{x}\left( t\right) =\left( 1+\lambda t\right) ^{x/\lambda
}=\sum\limits_{n=0}^{\infty }\left( x\right) _{n,\lambda }\frac{t^{n}}{n!}%
\text{ and }e_{\lambda }\left( t\right) =e_{\lambda }^{1}\left( t\right)
=\left( 1+\lambda t\right) ^{1/\lambda }\text{,}  \label{4}
\end{equation}%
where $\left( x\right) _{0,1}=1$ and $\left( x\right) _{n,\lambda }=x\left(
x-\lambda \right) \left( x-2\lambda \right) \cdots \left( x-\left(
n-1\right) \lambda \right) $, $n\geq 1$.

For $x\in 
\mathbb{R}
$ and $k$ nonnegative integer, the degenerate $\lambda $-Stirling
polynomials of the second kind are defined by \cite{5}%
\begin{equation}
\frac{\left( e_{\lambda }\left( t\right) -1\right) ^{k}}{k!}e_{\lambda
}^{x}\left( t\right) =\sum\limits_{n=k}^{\infty }S_{2,\lambda }^{\left(
x\right) }\left( n,k\right) \frac{t^{n}}{n!}\text{.}  \label{5}
\end{equation}

Note that%
\begin{equation*}
\underset{\lambda \longrightarrow 0}{\lim }\sum\limits_{n=k}^{\infty
}S_{2,\lambda }^{\left( x\right) }\left( n,k\right) \frac{t^{n}}{n!}=\frac{%
\left( e^{t}-1\right) ^{k}}{k!}e^{xt}\text{.}
\end{equation*}

From (\ref{4}), we get%
\begin{equation}
\left( t+x\right) _{n,\lambda }=\sum\limits_{k=0}^{n}S_{2,\lambda }^{\left(
x\right) }\left( n,k\right) \left( t\right) _{k}\text{, }n>0\text{,}
\label{6}
\end{equation}%
where $\left( t\right) _{0}=1$, $\left( t\right) _{n}=t\left( t-1\right)
\left( t-2\right) \cdots \left( t-\left( n-1\right) \right) $, $n\geq 1$.

Using (\ref{4}) and (\ref{6}), we note that%
\begin{equation*}
e_{\lambda }^{\left( x+y\right) }\left( t\right) =\sum\limits_{n=0}^{\infty
}\left( x+y\right) _{n,\lambda }\frac{t^{n}}{n!}=\sum\limits_{n=0}^{\infty
}\sum\limits_{k=0}^{n}S_{2,\lambda }^{\left( x\right) }\left( n,k\right)
\left( y\right) _{k}\frac{t^{n}}{n!}\text{.}
\end{equation*}

The degenerate Stirling numbers of the first kind are defined by (\cite{3}-%
\cite{10})%
\begin{equation}
\frac{1}{k!}\left( \log _{\lambda }\left( 1+t\right) \right)
^{k}=\sum\limits_{n=k}^{\infty }S_{1,\lambda }\left( n,k\right) \frac{t^{n}}{%
n!}\text{, }k\geq 0\text{.}  \label{7}
\end{equation}

Note here that $\underset{\lambda \longrightarrow 0}{\lim }S_{1,\lambda
}\left( n,l\right) =S_{1}\left( n,l\right) $ where $S_{1}\left( n,l\right) $
are the Stirling numbers of the first kind given by \cite{5}%
\begin{equation}
\frac{\left( \log \left( 1+t\right) \right) ^{k}}{k!}=\sum\limits_{n=k}^{%
\infty }S_{1}\left( n,k\right) \frac{t^{n}}{n!}\text{, }k\geq 0\text{.}
\label{8}
\end{equation}

The degenerate Stirling numbers of the second kind are defined by (\cite{3}-%
\cite{10})%
\begin{equation}
\frac{\left( e_{\lambda }\left( t\right) -1\right) ^{k}}{k!}%
=\sum\limits_{n=k}^{\infty }S_{2}\left( n,k\right) \frac{t^{n}}{n!}\text{, }%
k\geq 0\text{.}  \label{9}
\end{equation}

Observe that $\underset{\lambda \longrightarrow 0}{\lim }S_{2,\lambda
}\left( n,l\right) =S_{2}\left( n,l\right) $ where $S_{2}\left( n,l\right) $
are the Stirling numbers of the second kind given by \cite{5}%
\begin{equation}
\frac{\left( e^{t}-1\right) ^{k}}{k!}=\sum\limits_{n=k}^{\infty }S_{2}\left(
n,k\right) \frac{t^{n}}{n!}\text{, }k\geq 0\text{.}  \label{10}
\end{equation}

The degenerate Bernoulli polynomials of the second kind are given by (\cite%
{6}, \cite{8})%
\begin{equation}
\frac{t}{\log _{\lambda }\left( 1+t\right) }\left( 1+t\right)
^{x}=\sum\limits_{n=0}^{\infty }b_{n,\lambda }\left( x\right) \frac{t^{n}}{n!%
}\text{.}  \label{11}
\end{equation}

Note that $\underset{\lambda \longrightarrow 0}{\lim }b_{n,\lambda }\left(
x\right) =b_{n}\left( x\right) $ where $b_{n}\left( x\right) $ are the
Bernoulli polynomials of the second kind given by \cite{6}%
\begin{equation}
\frac{t}{\log \left( 1+t\right) }\left( 1+t\right)
^{x}=\sum\limits_{n=0}^{\infty }b_{n}\left( x\right) \frac{t^{n}}{n!}\text{.}
\label{12}
\end{equation}

\section{Degenerate Poly-Frobenius-Genocchi Numbers and Polynomials}

In this section, we introduce and investigate the modified polyexponential
functions. We give some identities and explicit relations for the modified
degenerate polyexponential functions. We define the degenerate
poly-Frobenius-Genocchi polynomials. Also, we give some relations and
identities for these polynomials.

In \cite{2}, Boyadzhiev introduced the polyexponential function, Kim et al.
in (\cite{6}, \cite{7}) considered and investigated the polyexponential
functions and the degenerate polyexponential functions.

The polyexponential functions are defined by (\cite{3}-\cite{11}, \cite{14})%
\begin{equation}
\func{Ei}_{k}\left( x\right) =\sum\limits_{n=1}^{\infty }\frac{x^{n}}{%
n^{k}\left( n-1\right) !}\text{, }k\in 
\mathbb{Z}
\text{.}  \label{13}
\end{equation}

For $k=1$, $\func{Ei}_{1}\left( x\right) =e^{x}-1$.

The modified degenerate polyexponential functions are given by (\cite{3}-%
\cite{11}, \cite{14})%
\begin{equation}
\func{Ei}_{k,\lambda }\left( x\right) =\sum\limits_{n=1}^{\infty }\frac{%
\left( 1\right) _{n,\lambda }}{n^{k}\left( n-1\right) !}x^{n}\text{, }%
\lambda \in 
\mathbb{R}
\text{.}  \label{14}
\end{equation}

Note that 
\begin{equation*}
\func{Ei}_{1,\lambda }\left( x\right) =\sum\limits_{n=1}^{\infty }\left(
1\right) _{n,\lambda }\frac{x^{n}}{n!}=e_{\lambda }\left( x\right) -1\text{.}
\end{equation*}

For $k\in 
\mathbb{Z}
$ and by means of the modified degenerate polyexponential functions. We
define the degenerate poly-Frobenius-Genocchi polynomials by the following
generating functions.%
\begin{equation}
\sum\limits_{n=0}^{\infty }FG_{n,\lambda }^{\left( k\right) }\left(
x,u\right) \frac{t^{n}}{n!}=\frac{\left( 1-u\right) \func{Ei}_{_{k,\lambda
}}\left( \log _{\lambda }\left( 1+t\right) \right) }{e_{\lambda }\left(
t\right) -u}e_{\lambda }^{x}\left( t\right) \text{.}  \label{15}
\end{equation}

When $x=0$, $FG_{n,\lambda }^{\left( k\right) }\left( u\right)
:=FG_{n,\lambda }^{\left( k\right) }\left( 0,u\right) $ are called the
degenerate poly-Frobenius-Genocchi numbers, where $log_{\lambda }\left(
t\right) =\frac{1}{\lambda }\left( t^{\lambda }-1\right) $ is the
compositional inverse of $e_{\lambda }\left( t\right) $ satisfying%
\begin{equation*}
\log _{\lambda }\left( e_{\lambda }\left( t\right) \right) =e_{\lambda
}\left( \log _{\lambda }\left( 1+t\right) \right) =t\text{.}
\end{equation*}

For $k=1$ and $u=-1$, we get the degenerate Genocchi polynomials%
\begin{equation*}
\sum\limits_{n=0}^{\infty }FG_{n,\lambda }^{\left( 1\right) }\left(
x,-1\right) \frac{t^{n}}{n!}=\frac{2\func{Ei}_{_{1,\lambda }}\left( \log
_{\lambda }\left( 1+t\right) \right) }{e_{\lambda }\left( t\right) +1}%
e_{\lambda }^{x}\left( t\right)
\end{equation*}%
\begin{equation*}
=\frac{2t}{e_{\lambda }\left( t\right) +1}e_{\lambda }^{x}\left( t\right)
=\sum\limits_{n=0}^{\infty }G_{n,\lambda }\left( x\right) \frac{t^{n}}{n!}%
\text{.}
\end{equation*}

From (\ref{15}), we can write the following equations%
\begin{equation}
FG_{n,\lambda }^{\left( k\right) }\left( x,u\right) =\sum\limits_{m=0}^{n}%
\binom{n}{m}FG_{n,\lambda }^{\left( k\right) }\text{ }\left( x\right)
_{n-m,\lambda }\text{,}  \tag{i}
\end{equation}%
\begin{eqnarray}
FG_{n,\lambda }^{\left( k\right) }\left( x+y,u\right)
&=&\sum\limits_{m=0}^{n}\binom{n}{m}FG_{n,\lambda }^{\left( k\right) }\left(
x,u\right) \left( y\right) _{n-m,\lambda }\text{,}  \TCItag{ii} \\
&=&\sum\limits_{m=0}^{n}\binom{n}{m}FG_{n,\lambda }^{\left( k\right) }\left(
y,u\right) \left( x\right) _{n-m,\lambda }\text{.}  \notag
\end{eqnarray}%
\begin{equation}
FG_{n,\lambda }^{\left( k\right) }\left( x+y,u\right) =\sum\limits_{m=0}^{n}%
\binom{n}{m}FG_{n,\lambda }^{\left( k\right) }\text{ }\left( x+y\right)
_{n-m,\lambda }\text{.}  \tag{iii}
\end{equation}

By (\ref{8}) and (\ref{14}), we get%
\begin{eqnarray}
\func{Ei}_{k,\lambda }\left( \log _{\lambda }\left( 1+t\right) \right)
&=&\sum\limits_{n=1}^{\infty }\frac{\left( 1\right) _{n,\lambda }\left( \log
_{\lambda }\left( 1+t\right) \right) ^{n}}{n^{k}\left( n-1\right) !}  \notag
\\
&=&\sum\limits_{n=1}^{\infty }\frac{\left( 1\right) _{n,\lambda }}{\left(
n-1\right) ^{k}}\sum\limits_{m=n}^{\infty }S_{1,\lambda }\left( m,n\right) 
\frac{t^{n}}{n!}  \notag \\
&=&t\sum\limits_{m=0}^{\infty }\sum\limits_{n=1}^{m+1}\frac{\left( 1\right)
_{n,\lambda }}{n^{k}-1}\frac{S_{1,\lambda }\left( m+1,n\right) }{m+1}\frac{%
t^{m}}{m!}\text{.}  \label{16}
\end{eqnarray}

Using (\ref{15}) and (\ref{16}), we get%
\begin{equation*}
\sum\limits_{n=0}^{\infty }FG_{n,\lambda }^{\left( k\right) }\left(
x,u\right) \frac{t^{n}}{n!}=\frac{\left( 1-u\right) e_{q}^{x}\left( t\right) 
}{e_{\lambda }\left( t\right) -u}\func{Ei}_{k,\lambda }\left( \log _{\lambda
}\left( 1+t\right) \right)
\end{equation*}%
\begin{equation*}
=t\sum\limits_{l=0}^{\infty }FG_{l,\lambda }\left( x,u\right) \frac{t^{l}}{l!%
}\sum\limits_{m=0}^{\infty }\frac{1}{m+1}\sum\limits_{j=1}^{m+1}\frac{\left(
1\right) _{j,\lambda }}{j^{k-1}}S_{1,\lambda }\left( m+1,j\right) \frac{t^{m}%
}{m!}\text{.}
\end{equation*}

By using Cauchy product and comparing the coefficients of $\frac{t^{n}}{n!}$
the above equations, we have the following theorem.

\begin{theorem}
For $n\geq 0$, we have%
\begin{equation*}
FG_{n,\lambda }^{\left( k\right) }\left( x,u\right) =n\sum\limits_{m=0}^{n-1}%
\frac{\binom{n-1}{m}}{m+1}\sum\limits_{j=1}^{m+1}\frac{\left( 1\right)
_{j,\lambda }}{j^{k-1}}S_{1,\lambda }\left( m+1,j\right) \text{ }%
FG_{n-1-m,\lambda }\left( x,u\right) \text{.}
\end{equation*}
\end{theorem}

From (\ref{15}), we write as%
\begin{equation*}
\sum\limits_{n=0}^{\infty }FG_{n,\lambda }^{\left( k\right) }\left(
x,u\right) \frac{t^{n}}{n!}\left( e_{q}\left( t\right) -u\right) =\left(
1-u\right) \func{Ei}_{k,\lambda }\left( \log _{\lambda }\left( 1+t\right)
\right) e_{\lambda }^{x}\left( t\right)
\end{equation*}%
and%
\begin{equation}
=\left( 1-u\right) \sum\limits_{n=0}^{\infty }\left( \sum\limits_{m=0}^{n}%
\binom{n}{m}\sum\limits_{j=1}^{m+1}\frac{\left( 1\right) _{j,\lambda }}{%
j^{k-1}}\frac{S_{1,\lambda }\left( m+1,j\right) }{m+1}\left( x\right)
_{n-m,\lambda }\right) \frac{t^{n+1}}{n!}\text{.}  \label{17}
\end{equation}

Comparing the coefficients of both sides in (\ref{17}).

We have the following theorem.

\begin{theorem}
For $n\geq 0$, we have%
\begin{eqnarray*}
&&\sum\limits_{m=0}^{n}\binom{n}{m}FG_{m,\lambda }^{\left( k\right) }\left(
x,u\right) \left( 1\right) _{n-m,\lambda }-u\text{ }FG_{n,\lambda }^{\left(
k\right) }\left( x,u\right) \\
&=&\left( 1-u\right) n\sum\limits_{m=0}^{n-1}\frac{\binom{n-1}{m}}{m+1}%
\sum\limits_{j=1}^{m+1}\frac{\left( 1\right) _{j,\lambda }}{j^{k-1}}%
S_{1,\lambda }\left( m+1,j\right) \left( x\right) _{n-1-m,\lambda }\text{.}
\end{eqnarray*}
\end{theorem}

From (\ref{14}), we note that%
\begin{equation}
\frac{d}{dx}\func{Ei}_{k,\lambda }\left( x\right) =\frac{d}{dx}%
\sum\limits_{n=1}^{\infty }\frac{\left( 1\right) _{j,\lambda }}{\left(
n-1\right) !n^{k}}x^{n}=\frac{1}{x}\func{Ei}_{k-1,\lambda }\left( x\right) 
\text{.}  \label{18}
\end{equation}

Thus, by (\ref{17}), we get%
\begin{equation*}
\func{Ei}_{k,\lambda }\left( x\right) =\int\limits_{0}^{x}\frac{1}{t}\func{Ei%
}_{k-1,\lambda }\left( t\right) dt
\end{equation*}%
\begin{equation*}
=\underset{\left( k-2\right) \text{ times}}{\underbrace{\int\limits_{0}^{x}%
\frac{1}{t}\int\limits_{0}^{t}\cdots \frac{1}{t}\int\limits_{0}^{t}\frac{1}{t%
}}}\func{Ei}_{1,\lambda }\left( x\right) \underset{\left( k-2\right) \text{
times}}{\underbrace{dt\cdots dt}}
\end{equation*}%
\begin{equation*}
=\underset{\left( k-2\right) \text{ times}}{\underbrace{\int\limits_{0}^{x}%
\frac{1}{t}\int\limits_{0}^{t}\cdots \frac{1}{t}\int\limits_{0}^{t}\frac{1}{t%
}}}\left( e_{\lambda }\left( t\right) -1\right) \underset{\left( k-2\right) 
\text{ times}}{\underbrace{dt\cdots dt}}\text{,}
\end{equation*}%
where $k\in 
\mathbb{Z}
^{+}$ with $k\geq 2$.

Form (\ref{11}), (\ref{15}) and (\ref{18}), for $k=2$%
\begin{equation*}
\sum\limits_{n=0}^{\infty }FG_{n,\lambda }^{\left( 2\right) }\left(
x,u\right) \frac{t^{n}}{n!}=\frac{1-u}{e_{\lambda }\left( t\right) -u}%
\int\limits_{0}^{t}\frac{t}{\log _{\lambda }\left( 1+t\right) }\left(
1+t\right) ^{\lambda -1}dt
\end{equation*}%
\begin{equation*}
=\frac{1-u}{e_{\lambda }\left( t\right) -u}\sum\limits_{m=0}^{\infty }\frac{%
b_{m,\lambda }\left( \lambda -1\right) }{m+1}\frac{t^{m}}{m!}%
=\sum\limits_{l=0}^{\infty }H_{l,\lambda }\left( 0,u\right) \frac{t^{l}}{l!}%
\sum\limits_{m=0}^{\infty }\frac{b_{m,\lambda }\left( \lambda -1\right) }{m+1%
}\frac{t^{m}}{m!}\text{.}
\end{equation*}%
From the last equations, we have the following theorem.

\begin{theorem}
For $n\geq 0$, we have%
\begin{equation*}
FG_{n,\lambda }^{\left( 2\right) }\left( x,u\right) =\sum\limits_{m=0}^{n}%
\binom{n}{m}H_{n-m,\lambda }\left( 0,u\right) \frac{b_{m,\lambda }\left(
\lambda -1\right) }{m+1}\text{,}
\end{equation*}%
where $H_{n,\lambda }\left( 0,u\right) $ is degenerate Frobenius-Euler
numbers.
\end{theorem}

Recently, Masjed-Jamai et al in \cite{13} and Srivastava et al in (\cite{15}%
, \cite{16}) introduced a new type parametric Euler numbers and polynomials
as%
\begin{equation*}
\frac{2}{e^{t}+1}e^{pt}\cos \left( qt\right) =\sum\limits_{n=0}^{\infty
}E_{n}^{\left( c\right) }\left( p,q\right) \frac{t^{n}}{n!}
\end{equation*}%
and%
\begin{equation*}
\frac{2}{e^{t}+1}e^{pt}\sin \left( qt\right) =\sum\limits_{n=0}^{\infty
}E_{n}^{\left( s\right) }\left( p,q\right) \frac{t^{n}}{n!}\text{,}
\end{equation*}%
where%
\begin{equation*}
e^{pt}\cos \left( qt\right) =\sum\limits_{n=0}^{\infty }C_{n}\left(
p,q\right) \frac{t^{n}}{n!}
\end{equation*}%
and%
\begin{equation*}
e^{pt}\sin \left( qt\right) =\sum\limits_{n=0}^{\infty }S_{n}\left(
p,q\right) \frac{t^{n}}{n!}\text{.}
\end{equation*}

\section{Degenerate Poly-Frobenius-Genocchi Polynomials of Complex Variables}

In this section, we define the Frobenius-Genocchi polynomials of the complex
variables. We consider the degenerate cosine function and the degenerate
sine function. Using the degenerate cosine function and the degenerate sine
function, we introduce the cosine degenerate poly-Frobenius-Genocchi
polynomials and the sine degenerate poly-Frobenius-Genocchi polynomials.

From (\ref{15}), we write as%
\begin{equation*}
\sum\limits_{n=0}^{\infty }FG_{n,\lambda }^{\left( k\right) }\left(
x+iy;u\right) \frac{t^{n}}{n!}=\frac{\left( 1-u\right) \func{Ei}%
_{_{k,\lambda }}\left( \log _{\lambda }\left( 1+t\right) \right) }{%
e_{q}\left( t\right) -u}e_{\lambda }^{\left( x+iy\right) }\left( t\right)
\end{equation*}%
\begin{equation}
=\frac{\left( 1-u\right) \func{Ei}_{_{k,\lambda }}\left( \log _{\lambda
}\left( 1+t\right) \right) }{e_{q}\left( t\right) -u}e_{\lambda }^{x}\left(
t\right) \left[ \cos _{\lambda }^{\left( y\right) }\left( t\right) +i\sin
_{\lambda }^{\left( y\right) }\left( t\right) \right]  \label{19}
\end{equation}%
and%
\begin{eqnarray}
&&\sum\limits_{n=0}^{\infty }FG_{n,\lambda }^{\left( k\right) }\left(
x-iy;u\right) \frac{t^{n}}{n!}  \notag \\
&=&\frac{\left( 1-u\right) \func{Ei}_{_{k,\lambda }}\left( \log _{\lambda
}\left( 1+t\right) \right) }{e_{q}\left( t\right) -u}e_{\lambda }^{x}\left(
t\right) \left[ \cos _{\lambda }^{\left( y\right) }\left( t\right) -i\sin
_{\lambda }^{\left( y\right) }\left( t\right) \right] \text{.}  \label{20}
\end{eqnarray}

By (\ref{19}) and (\ref{20}), we get%
\begin{eqnarray}
&&\frac{\left( 1-u\right) \func{Ei}_{_{k,\lambda }}\left( \log _{\lambda
}\left( 1+t\right) \right) }{e_{q}\left( t\right) -u}e_{\lambda }^{x}\left(
t\right) \cos _{\lambda }^{\left( y\right) }\left( t\right)  \notag \\
&=&\sum\limits_{n=0}^{\infty }\frac{FG_{n,\lambda }^{\left( k\right) }\left(
x+iy;u\right) +FG_{n,\lambda }^{\left( k\right) }\left( x-iy;u\right) }{2}%
\frac{t^{n}}{n!}  \label{21}
\end{eqnarray}%
and%
\begin{eqnarray}
&&\frac{\left( 1-u\right) \func{Ei}_{_{k,\lambda }}\left( \log _{\lambda
}\left( 1+t\right) \right) }{e_{q}\left( t\right) -u}e_{\lambda }^{x}\left(
t\right) \sin _{\lambda }^{\left( y\right) }\left( t\right)  \notag \\
&=&\sum\limits_{n=0}^{\infty }\frac{FG_{n,\lambda }^{\left( k\right) }\left(
x+iy;u\right) -FG_{n,\lambda }^{\left( k\right) }\left( x-iy;u\right) }{2i}%
\frac{t^{n}}{n!}\text{.}  \label{22}
\end{eqnarray}

Using (\ref{4}), we define the degenerate cosine-functions and the
degenerate sine-functions as%
\begin{equation}
\cos _{\lambda }^{\left( y\right) }\left( t\right) =\frac{e_{\lambda
}^{\left( iy\right) }\left( t\right) +e_{\lambda }^{\left( -iy\right)
}\left( t\right) }{2}=\cos \left( \frac{y}{\lambda }\log \left( 1+\lambda
t\right) \right)  \label{23}
\end{equation}%
and%
\begin{equation}
\sin _{\lambda }^{\left( y\right) }\left( t\right) =\frac{e_{\lambda
}^{\left( iy\right) }\left( t\right) -e_{\lambda }^{\left( -iy\right)
}\left( t\right) }{2i}=\sin \left( \frac{y}{\lambda }\log \left( 1+\lambda
t\right) \right) \text{,}  \label{24}
\end{equation}%
where $\underset{\lambda \longrightarrow 0}{\lim }\cos _{\lambda }^{\left(
y\right) }\left( t\right) =\cos (yt)$ and $\underset{\lambda \longrightarrow
0}{\lim }\sin _{\lambda }^{\left( y\right) }\left( t\right) =\sin (yt)$.

Now, we define the cosine degenerate poly-Frobenius-Genocchi polynomials and
the sine degenerate poly-Frobenius-Genocchi polynomials, respectively;%
\begin{equation}
\sum\limits_{n=0}^{\infty }FG_{n,\lambda }^{\left[ k,c\right] }\left(
x,y;u\right) \frac{t^{n}}{n!}=\frac{\left( 1-u\right) \func{Ei}_{_{k,\lambda
}}\left( \log _{\lambda }\left( 1+t\right) \right) }{e_{\lambda }\left(
t\right) -u}e_{\lambda }^{x}\left( t\right) \cos _{\lambda }^{\left(
y\right) }\left( t\right)  \label{25}
\end{equation}%
and%
\begin{equation}
\sum\limits_{n=0}^{\infty }FG_{n,\lambda }^{\left[ k,s\right] }\left(
x,y;u\right) \frac{t^{n}}{n!}=\frac{\left( 1-u\right) \func{Ei}_{_{k,\lambda
}}\left( \log _{\lambda }\left( 1+t\right) \right) }{e_{\lambda }\left(
t\right) -u}e_{\lambda }^{x}\left( t\right) \sin _{\lambda }^{\left(
y\right) }\left( t\right) \text{.}  \label{26}
\end{equation}

From (\ref{4}), we write%
\begin{equation*}
e_{\lambda }^{\left( iy\right) }\left( t\right) =\sum\limits_{n=0}^{\infty
}\left( iy\right) _{n,\lambda }\frac{t^{n}}{n!}\text{ and }e_{\lambda
}^{\left( -iy\right) }\left( t\right) =\sum\limits_{n=0}^{\infty }\left(
-iy\right) _{n,\lambda }\frac{t^{n}}{n!}\text{.}
\end{equation*}

Using (\ref{23}) and (\ref{24}), we get%
\begin{equation}
\cos _{\lambda }^{\left( y\right) }\left( t\right) =\frac{1}{2}%
\sum\limits_{n=0}^{\infty }\left( \left( iy\right) _{n,\lambda }+\left(
-iy\right) _{n,\lambda }\right) \frac{t^{n}}{n!}  \label{27}
\end{equation}%
and%
\begin{equation}
\sin _{\lambda }^{\left( y\right) }\left( t\right) =\frac{1}{2i}%
\sum\limits_{n=0}^{\infty }\left( \left( iy\right) _{n,\lambda }-\left(
-iy\right) _{n,\lambda }\right) \frac{t^{n}}{n!}\text{.}  \label{28}
\end{equation}

By (\ref{4}), (\ref{27}) and (\ref{4}), (\ref{28}), we have the following
equations, respectively,%
\begin{equation}
e_{\lambda }^{x}\left( t\right) \cos _{\lambda }^{\left( y\right) }\left(
t\right) =\frac{1}{2}\sum\limits_{n=0}^{\infty }\sum\limits_{k=0}^{n}\binom{n%
}{k}\left( x\right) _{n-k,\lambda }\left( \left( iy\right) _{n,\lambda
}+\left( -iy\right) _{n,\lambda }\right) \frac{t^{n}}{n!}  \label{29}
\end{equation}%
and%
\begin{equation}
e_{\lambda }^{x}\left( t\right) \sin _{\lambda }^{\left( y\right) }\left(
t\right) =\frac{1}{2i}\sum\limits_{n=0}^{\infty }\sum\limits_{k=0}^{n}\binom{%
n}{k}\left( x\right) _{n-k,\lambda }\left( \left( iy\right) _{n,\lambda
}-\left( -iy\right) _{n,\lambda }\right) \frac{t^{n}}{n!}\text{.}  \label{30}
\end{equation}

From (\ref{25}) and (\ref{29}), we write%
\begin{equation*}
\sum\limits_{n=0}^{\infty }FG_{n,\lambda }^{\left[ k,c\right] }\left(
x,y;u\right) \frac{t^{n}}{n!}=\frac{\left( 1-u\right) \func{Ei}_{_{k,\lambda
}}\left( \log _{\lambda }\left( 1+t\right) \right) }{e_{q}\left( t\right) -u}%
e_{q}^{x}\left( t\right) \cos _{\lambda }^{\left( y\right) }\left( t\right)
\end{equation*}%
\begin{equation*}
=\sum\limits_{n=0}^{\infty }FG_{n,\lambda }^{\left( k\right) }\frac{t^{n}}{n!%
}\frac{1}{2}\sum\limits_{n=0}^{\infty }\sum\limits_{k=0}^{n}\binom{n}{k}%
\left( x\right) _{n-k,\lambda }\left( \left( iy\right) _{n,\lambda }+\left(
-iy\right) _{n,\lambda }\right) \frac{t^{n}}{n!}\text{.}
\end{equation*}

Using the Cauchy product and comparing the coefficients, we have%
\begin{equation}
FG_{n,\lambda }^{\left[ k,c\right] }\left( x,y;u\right) =\frac{1}{2}%
\sum\limits_{j=0}^{n}\binom{n}{j}\text{ }FG_{n-j,\lambda }^{\left( k\right)
}\sum\limits_{k=0}^{j}\binom{j}{k}\left( x\right) _{j-k,\lambda }\left(
\left( iy\right) _{k,\lambda }+\left( -iy\right) _{k,\lambda }\right) \text{.%
}  \label{32}
\end{equation}

From (\ref{26}) and (\ref{29}), similarly, we have%
\begin{equation}
FG_{n,\lambda }^{\left[ k,s\right] }\left( x,y;u\right) =\frac{1}{2i}%
\sum\limits_{j=0}^{n}\binom{n}{j}\text{ }FG_{n-j,\lambda }^{\left( k\right)
}\sum\limits_{k=0}^{j}\binom{j}{k}\left( x\right) _{j-k,\lambda }\left(
\left( iy\right) _{k,\lambda }-\left( -iy\right) _{k,\lambda }\right) \text{.%
}  \label{33}
\end{equation}

From (\ref{32}) and (\ref{33}), we have the following theorems.

\begin{theorem}
The following relations hold true:%
\begin{equation*}
FG_{n,\lambda }^{\left[ k,c\right] }\left( x,y;u\right) =\frac{1}{2}%
\sum\limits_{j=0}^{n}\binom{n}{j}\text{ }FG_{n-j,\lambda }^{\left( k\right)
}\sum\limits_{k=0}^{j}\binom{j}{k}\left( x\right) _{j-k,\lambda }\left(
\left( iy\right) _{k,\lambda }+\left( -iy\right) _{k,\lambda }\right)
\end{equation*}%
and%
\begin{equation*}
FG_{n,\lambda }^{\left[ k,s\right] }\left( x,y;u\right) =\frac{1}{2i}%
\sum\limits_{j=0}^{n}\binom{n}{j}\text{ }FG_{n-j,\lambda }^{\left( k\right)
}\sum\limits_{k=0}^{j}\binom{j}{k}\left( x\right) _{j-k,\lambda }\left(
\left( iy\right) _{k,\lambda }-\left( -iy\right) _{k,\lambda }\right) \text{.%
}
\end{equation*}
\end{theorem}

Now, we define the degenerate two parametric $C_{n,\lambda }\left(
x,y\right) $ and $S_{n,\lambda }\left( x,y\right) $ polynomials,
respectively,%
\begin{equation}
e_{\lambda }^{\left( x\right) }\left( t\right) \cos _{\lambda }^{\left(
y\right) }\left( t\right) =\sum\limits_{n=0}^{\infty }C_{n,\lambda }\left(
x,y\right) \frac{t^{n}}{n!}  \label{34}
\end{equation}%
and%
\begin{equation}
e_{\lambda }^{\left( x\right) }\left( t\right) \sin _{\lambda }^{\left(
y\right) }\left( t\right) =\sum\limits_{n=0}^{\infty }S_{n,\lambda }\left(
x,y\right) \frac{t^{n}}{n!}\text{.}  \label{35}
\end{equation}

From (\ref{4}) and (\ref{27}), we get%
\begin{equation*}
C_{n,\lambda }\left( x,y\right) =\frac{1}{2}\sum\limits_{k=0}^{n}\binom{n}{k}%
\left( x\right) _{n-k,\lambda }\left( \left( iy\right) _{k,\lambda }+\left(
-iy\right) _{k,\lambda }\right) \text{.}
\end{equation*}

Similarly, (\ref{4}) and (\ref{28}), we get%
\begin{equation*}
S_{n,\lambda }\left( x,y\right) =\frac{1}{2i}\sum\limits_{k=0}^{n}\binom{n}{k%
}\left( x\right) _{n-k,\lambda }\left( \left( iy\right) _{k,\lambda }-\left(
-iy\right) _{k,\lambda }\right) \text{.}
\end{equation*}

From (\ref{16}), (\ref{25}) and (\ref{29}), we write%
\begin{eqnarray*}
&&\left( e_{q}\left( t\right) -u\right) \sum\limits_{n=0}^{\infty
}FG_{n,\lambda }^{\left[ k,c\right] }\left( x,y;u\right) \frac{t^{n}}{n!} \\
&=&\left( 1-u\right) \func{Ei}_{k,\lambda }\left( \log _{\lambda }\left(
1+t\right) \right) e_{q}^{\left( x\right) }\left( t\right) \cos _{\lambda
}^{\left( y\right) }\left( t\right) \text{.}
\end{eqnarray*}

The left hand side of this equation is%
\begin{equation}
\sum\limits_{n=0}^{\infty }\left\{ \sum\limits_{l=0}^{n}\binom{n}{l}\left(
1\right) _{n-l,\lambda }\text{ }FG_{l,\lambda }^{\left[ k,c\right] }\left(
x,y;u\right) -u\text{ }FG_{n,\lambda }^{\left[ k,c\right] }\left(
x,y;u\right) \right\} \frac{t^{n}}{n!}\text{.}  \label{36}
\end{equation}

The right hand side of this equation is%
\begin{eqnarray}
&&\frac{1-u}{2}\sum\limits_{n=0}^{\infty }n\sum\limits_{m=0}^{n-1}\frac{%
\binom{n-1}{m}}{m+1}\sum\limits_{j=0}^{m+1}\frac{\left( 1\right) _{j,\lambda
}}{j^{k-1}}S_{1,\lambda }\left( m+1,j\right)  \notag \\
&&\times \sum\limits_{l=0}^{n-1-m}\binom{n-1-m}{l}\left( x\right)
_{n-1-m,\lambda }\left( \left( iy\right) _{n-1-m,\lambda }+\left( -iy\right)
_{n-1-m,\lambda }\right) \frac{t^{n}}{n!}\text{.}  \label{37}
\end{eqnarray}

From (\ref{36}) and (\ref{37}), we get%
\begin{equation*}
2\left\{ \sum\limits_{l=0}^{n}\binom{n}{l}\left( 1\right) _{n-l,\lambda }%
\text{ }FG_{l,\lambda }^{\left[ k,c\right] }\left( x,y;u\right) -u\text{ }%
FG_{n,\lambda }^{\left[ k,c\right] }\left( x,y;u\right) \right\}
\end{equation*}%
\begin{eqnarray}
&=&\left( 1-u\right) \left\{ n\sum\limits_{m=0}^{n-1}\frac{\binom{n-1}{m}}{%
m+1}\sum\limits_{j=0}^{m+1}\frac{\left( 1\right) _{j,\lambda }}{j^{k-1}}%
S_{1,\lambda }\left( m+1,j\right) \right.  \notag \\
&&\times \left. \sum\limits_{l=0}^{n-1-m}\binom{n-1-m}{l}\left( x\right)
_{n-1-m,\lambda }\left( \left( iy\right) _{n-1-m,\lambda }+\left( -iy\right)
_{n-1-m,\lambda }\right) \right\} \text{.}  \label{38}
\end{eqnarray}

Similarly, (\ref{16}), (\ref{26}) and (\ref{30})%
\begin{equation*}
2i\left\{ \sum\limits_{l=0}^{n}\binom{n}{l}\left( 1\right) _{n-l,\lambda }%
\text{ }FG_{l,\lambda }^{\left[ k,s\right] }\left( x,y;u\right) -u\text{ }%
FG_{n,\lambda }^{\left[ k,s\right] }\left( x,y;u\right) \right\}
\end{equation*}%
\begin{eqnarray}
&=&\left( 1-u\right) \left\{ n\sum\limits_{m=0}^{n-1}\frac{\binom{n-1}{m}}{%
m+1}\sum\limits_{j=0}^{m+1}\frac{\left( 1\right) _{j,\lambda }}{j^{k-1}}%
S_{1,\lambda }\left( m+1,j\right) \right.  \notag \\
&&\times \left. \sum\limits_{l=0}^{n-1-m}\binom{n-1-m}{l}\left( x\right)
_{n-1-m,\lambda }\left( \left( iy\right) _{n-1-m,\lambda }-\left( -iy\right)
_{n-1-m,\lambda }\right) \right\} \text{.}  \label{39}
\end{eqnarray}

\begin{theorem}
The following relations hold true:%
\begin{equation*}
2\left\{ \sum\limits_{l=0}^{n}\binom{n}{l}\left( 1\right) _{n-l,\lambda }%
\text{ }FG_{l,\lambda }^{\left[ k,c\right] }\left( x,y;u\right) -u\text{ }%
FG_{n,\lambda }^{\left[ k,c\right] }\left( x,y;u\right) \right\}
\end{equation*}%
\begin{eqnarray*}
&=&\left( 1-u\right) \left\{ n\sum\limits_{m=0}^{n-1}\frac{\binom{n-1}{m}}{%
m+1}\sum\limits_{j=0}^{m+1}\frac{\left( 1\right) _{j,\lambda }}{j^{k-1}}%
S_{1,\lambda }\left( m+1,j\right) \right. \\
&&\times \left. \sum\limits_{l=0}^{n-1-m}\binom{n-1-m}{l}\left( x\right)
_{n-1-m,\lambda }\left( \left( iy\right) _{n-1-m,\lambda }+\left( -iy\right)
_{n-1-m,\lambda }\right) \right\}
\end{eqnarray*}%
and%
\begin{equation*}
2i\left\{ \sum\limits_{l=0}^{n}\binom{n}{l}\left( 1\right) _{n-l,\lambda }%
\text{ }FG_{l,\lambda }^{\left[ k,s\right] }\left( x,y;u\right) -u\text{ }%
FG_{n,\lambda }^{\left[ k,s\right] }\left( x,y;u\right) \right\}
\end{equation*}%
\begin{eqnarray*}
&=&\left( 1-u\right) \left\{ n\sum\limits_{m=0}^{n-1}\frac{\binom{n-1}{m}}{%
m+1}\sum\limits_{j=0}^{m+1}\frac{\left( 1\right) _{j,\lambda }}{j^{k-1}}%
S_{1,\lambda }\left( m+1,j\right) \right. \\
&&\times \left. \sum\limits_{l=0}^{n-1-m}\binom{n-1-m}{l}\left( x\right)
_{n-1-m,\lambda }\left( \left( iy\right) _{n-1-m,\lambda }-\left( -iy\right)
_{n-1-m,\lambda }\right) \right\} \text{.}
\end{eqnarray*}
\end{theorem}

From (\ref{6}) and (\ref{25}),%
\begin{equation*}
\sum\limits_{n=0}^{\infty }FG_{n,\lambda }^{\left[ k,c\right] }\left(
x_{1}+x_{2},y;u\right) \frac{t^{n}}{n!}=\frac{e_{\lambda }^{\left(
x_{1}+x_{2}\right) }\left( t\right) \left( 1-u\right) \func{Ei}_{k,\lambda
}\left( \log _{\lambda }\left( 1+t\right) \right) }{e_{q}\left( t\right) -u}%
\cos _{\lambda }^{\left( y\right) }\left( t\right)
\end{equation*}%
\begin{equation*}
=\sum\limits_{m=0}^{\infty }\sum\limits_{k=0}^{m}S_{2,\lambda }^{\left(
x_{1}\right) }\left( m,k\right) \left( x_{2}\right) _{k}\frac{t^{m}}{m!}%
\sum\limits_{l=0}^{\infty }FG_{l,\lambda }^{\left[ k,c\right] }\left(
0,y;u\right) \frac{t^{l}}{l!}
\end{equation*}%
\begin{equation}
=\sum\limits_{n=0}^{\infty }\left\{ \sum\limits_{m=0}^{n}\binom{n}{m}%
\sum\limits_{k=0}^{m}S_{2,\lambda }^{\left( x_{1}\right) }\left( m,k\right)
\left( x_{2}\right) _{k}\text{ }FG_{n-m,\lambda }^{\left[ k,c\right] }\left(
0,y;u\right) \right\} \frac{t^{n}}{n!}\text{.}  \label{40}
\end{equation}

From (\ref{6}) and (\ref{26}), we get%
\begin{equation*}
\sum\limits_{n=0}^{\infty }FG_{n,\lambda }^{\left[ k,s\right] }\left(
x_{1}+x_{2},y;u\right) \frac{t^{n}}{n!}=e_{\lambda }^{\left(
x_{1}+x_{2}\right) }\left( t\right) \frac{\left( 1-u\right) \func{Ei}%
_{k,\lambda }\left( \log _{\lambda }\left( 1+t\right) \right) }{e_{q}\left(
t\right) -u}\sin _{\lambda }^{\left( y\right) }\left( t\right)
\end{equation*}%
\begin{equation}
=\sum\limits_{n=0}^{\infty }\left\{ \sum\limits_{m=0}^{n}\binom{n}{m}%
\sum\limits_{k=0}^{m}S_{2,\lambda }^{\left( x_{1}\right) }\left( m,k\right)
\left( x_{2}\right) _{k}\text{ }FG_{n-m,\lambda }^{\left[ k,s\right] }\left(
0,y;u\right) \right\} \frac{t^{n}}{n!}\text{.}  \label{41}
\end{equation}

Comparing the coefficients of $\frac{t^{n}}{n!}$ both sides the equations (%
\ref{40}) and (\ref{41}), we have the following theorem.

\begin{theorem}
The following relations hold true:%
\begin{equation*}
FG_{n,\lambda }^{\left[ k,c\right] }\left( x_{1}+x_{2},y;u\right)
=\sum\limits_{m=0}^{n}\binom{n}{m}\sum\limits_{k=0}^{m}S_{2,\lambda
}^{\left( x_{1}\right) }\left( m,k\right) \left( x_{2}\right) _{k}\text{ }%
FG_{n-m,\lambda }^{\left[ k,c\right] }\left( 0,y;u\right)
\end{equation*}%
and%
\begin{equation*}
FG_{n,\lambda }^{\left[ k,s\right] }\left( x_{1}+x_{2},y;u\right)
=\sum\limits_{m=0}^{n}\binom{n}{m}\sum\limits_{k=0}^{m}S_{2,\lambda
}^{\left( x_{1}\right) }\left( m,k\right) \left( x_{2}\right) _{k}\text{ }%
FG_{n-m,\lambda }^{\left[ k,s\right] }\left( 0,y;u\right) \text{.}
\end{equation*}
\end{theorem}

Now, for our use in the present investigation, we find the expressions of $%
\cos _{\lambda }^{\left( x_{1}+x_{2}\right) }\left( t\right) $ and $\sin
_{\lambda }^{\left( x_{1}+x_{2}\right) }\left( t\right) $.

From (\ref{23}), we get%
\begin{equation*}
\cos _{\lambda }^{\left( x_{1}+x_{2}\right) }\left( t\right) =\cos \left( 
\frac{x_{1}+x_{2}}{\lambda }\log \left( 1+\lambda t\right) \right)
\end{equation*}%
\begin{eqnarray*}
&=&\cos \left( \frac{x_{1}}{\lambda }\log \left( 1+\lambda t\right) \right)
\cos \left( \frac{x_{2}}{\lambda }\log \left( 1+\lambda t\right) \right) \\
&&-\sin \left( \frac{x_{1}}{\lambda }\log \left( 1+\lambda t\right) \right)
\sin \left( \frac{x_{2}}{\lambda }\log \left( 1+\lambda t\right) \right)
\end{eqnarray*}%
\begin{equation}
=\cos _{\lambda }^{\left( x_{1}\right) }\left( t\right) \cos _{\lambda
}^{\left( x_{2}\right) }\left( t\right) -\sin _{\lambda }^{\left(
x_{1}\right) }\left( t\right) \sin _{\lambda }^{\left( x_{2}\right) }\left(
t\right) \text{.}  \label{41a}
\end{equation}

Putting (\ref{41a}), $x_{1}=x_{2}=x$, we get%
\begin{equation*}
\cos _{\lambda }^{\left( 2x\right) }\left( t\right) =\left( \cos _{\lambda
}^{\left( x\right) }\left( t\right) \right) ^{2}-\left( \sin _{\lambda
}^{\left( x\right) }\left( t\right) \right) ^{2}\text{.}
\end{equation*}

By (\ref{24}), we get%
\begin{equation*}
\sin _{\lambda }^{\left( x_{1}+x_{2}\right) }\left( t\right) =\sin \left( 
\frac{x_{1}+x_{2}}{\lambda }\log \left( 1+\lambda t\right) \right)
\end{equation*}%
\begin{equation}
=\sin _{\lambda }^{\left( x_{1}\right) }\left( t\right) \cos _{\lambda
}^{\left( x_{2}\right) }\left( t\right) +\cos _{\lambda }^{\left(
x_{1}\right) }\left( t\right) \sin _{\lambda }^{\left( x_{2}\right) }\left(
t\right) \text{.}  \label{42}
\end{equation}

Setting (\ref{42}), $x_{1}=x_{2}=x$, we get%
\begin{equation*}
\sin _{\lambda }^{\left( 2x\right) }\left( t\right) =2\cos _{\lambda
}^{\left( x\right) }\left( t\right) \sin _{\lambda }^{\left( x\right)
}\left( t\right) \text{.}
\end{equation*}

From (\ref{34}) and (\ref{41a}), we write%
\begin{equation*}
\sum\limits_{n=0}^{\infty }C_{n,\lambda }\left(
x_{1}+x_{2},y_{1}+y_{2}\right) \frac{t^{n}}{n!}=e_{\lambda }^{\left(
x_{1}+x_{2}\right) }\left( t\right) \cos _{\lambda }^{\left(
y_{1}+y_{2}\right) }\left( t\right)
\end{equation*}%
\begin{equation*}
=e_{\lambda }^{\left( x_{1}\right) }\left( t\right) e_{\lambda }^{\left(
x_{2}\right) }\left( t\right) \left( \cos _{\lambda }^{\left( y_{1}\right)
}\left( t\right) \cos _{\lambda }^{\left( y_{2}\right) }\left( t\right)
-\sin _{\lambda }^{\left( y_{1}\right) }\left( t\right) \sin _{\lambda
}^{\left( y_{2}\right) }\left( t\right) \right)
\end{equation*}%
\begin{equation}
=\sum\limits_{n=0}^{\infty }C_{n,\lambda }\left( x_{1},y_{1}\right) \frac{%
t^{n}}{n!}\sum\limits_{n=0}^{\infty }C_{n,\lambda }\left( x_{2},y_{2}\right) 
\frac{t^{n}}{n!}-\sum\limits_{n=0}^{\infty }S_{n,\lambda }\left(
x_{1},y_{1}\right) \frac{t^{n}}{n!}\sum\limits_{n=0}^{\infty }S_{n,\lambda
}\left( x_{2},y_{2}\right) \frac{t^{n}}{n!}\text{.}  \label{43}
\end{equation}

Using (\ref{35}) and (\ref{42}), we write%
\begin{equation*}
\sum\limits_{n=0}^{\infty }S_{n,\lambda }\left(
x_{1}+x_{2},y_{1}+y_{2}\right) \frac{t^{n}}{n!}=e_{\lambda }^{\left(
x_{1}+x_{2}\right) }\left( t\right) \sin _{\lambda }^{\left(
y_{1}+y_{2}\right) }\left( t\right)
\end{equation*}%
\begin{equation}
=\sum\limits_{n=0}^{\infty }S_{n,\lambda }\left( x_{1},y_{1}\right) \frac{%
t^{n}}{n!}\sum\limits_{n=0}^{\infty }C_{n,\lambda }\left( x_{2},y_{2}\right) 
\frac{t^{n}}{n!}+\sum\limits_{n=0}^{\infty }C_{n,\lambda }\left(
x_{1},y_{1}\right) \frac{t^{n}}{n!}\sum\limits_{n=0}^{\infty }S_{n,\lambda
}\left( x_{2},y_{2}\right) \frac{t^{n}}{n!}\text{.}  \label{44}
\end{equation}

By using Cauchy product above the equations (\ref{43}) and (\ref{44}), we
have the following theorem.

\begin{theorem}
The following relations hold true:%
\begin{equation}
C_{n,\lambda }\left( x_{1}+x_{2},y_{1}+y_{2}\right) =\sum\limits_{k=0}^{n}%
\binom{n}{k}\left\{ C_{n-k,\lambda }\left( x_{1},y_{1}\right) C_{k,\lambda
}\left( x_{2},y_{2}\right) -S_{n-k,\lambda }\left( x_{1},y_{1}\right)
S_{k,\lambda }\left( x_{2},y_{2}\right) \right\}  \label{45}
\end{equation}%
and%
\begin{equation}
S_{n,\lambda }\left( x_{1}+x_{2},y_{1}+y_{2}\right) =\sum\limits_{k=0}^{n}%
\binom{n}{k}\left\{ S_{n-k,\lambda }\left( x_{1},y_{1}\right) C_{k,\lambda
}\left( x_{2},y_{2}\right) +C_{n-k,\lambda }\left( x_{1},y_{1}\right)
S_{k,\lambda }\left( x_{2},y_{2}\right) \right\} \text{.}  \label{46}
\end{equation}
\end{theorem}

Setting $x_{1}=x_{2}=x$ and $y_{1}=y_{2}=y$ in (\ref{45}) and (\ref{46}), we
have respectively,%
\begin{equation*}
C_{n,\lambda }\left( 2x,2y\right) =\sum\limits_{k=0}^{n}\binom{n}{k}\left\{
C_{n-k,\lambda }\left( x,y\right) C_{k,\lambda }\left( x,y\right)
-S_{n-k,\lambda }\left( x,y\right) S_{k,\lambda }\left( x,y\right) \right\}
\end{equation*}%
and%
\begin{equation*}
S_{n,\lambda }\left( 2x,2y\right) =2\sum\limits_{k=0}^{n}\binom{n}{k}%
S_{n-k,\lambda }\left( x,y\right) C_{k,\lambda }\left( x,y\right) \text{.}
\end{equation*}

From (\ref{25}) and (\ref{41}), we write%
\begin{equation*}
\sum\limits_{n=0}^{\infty }FG_{n,\lambda }^{\left[ k,c\right] }\left(
x_{1}+x_{2},y_{1}+y_{2};u\right) \frac{t^{n}}{n!}
\end{equation*}%
\begin{eqnarray*}
&=&\frac{\left( 1-u\right) \func{Ei}_{k,\lambda }\left( \log _{\lambda
}\left( 1+t\right) \right) }{e_{q}\left( t\right) -u}e_{\lambda }^{\left(
x_{1}\right) }\left( t\right) e_{\lambda }^{\left( x_{2}\right) }\left(
t\right) \\
&&\left\{ \cos _{\lambda }^{\left( x_{1}\right) }\left( t\right) \cos
_{\lambda }^{\left( x_{2}\right) }\left( t\right) -\sin _{\lambda }^{\left(
y_{1}\right) }\left( t\right) \sin _{\lambda }^{\left( y_{2}\right) }\left(
t\right) \right\}
\end{eqnarray*}%
\begin{eqnarray}
&=&\sum\limits_{m=0}^{\infty }FG_{m,\lambda }^{\left[ k,c\right] }\left(
x_{1},y_{1};u\right) \frac{t^{m}}{m!}\sum\limits_{k=0}^{\infty }C_{k,\lambda
}\left( x_{2},y_{2}\right) \frac{t^{k}}{k!}  \notag \\
&&-\sum\limits_{m=0}^{\infty }FG_{m,\lambda }^{\left[ k,s\right] }\left(
x_{1},y_{1};u\right) \frac{t^{m}}{m!}\sum\limits_{k=0}^{\infty }S_{k,\lambda
}\left( x_{2},y_{2}\right) \frac{t^{k}}{k!}\text{.}  \label{47}
\end{eqnarray}

Using (\ref{26}) and (\ref{42}), we write%
\begin{equation*}
\sum\limits_{n=0}^{\infty }FG_{n,\lambda }^{\left[ k,s\right] }\left(
x_{1}+x_{2},y_{1}+y_{2};u\right) \frac{t^{n}}{n!}
\end{equation*}%
\begin{eqnarray*}
&=&\frac{\left( 1-u\right) \func{Ei}_{k,\lambda }\left( \log _{\lambda
}\left( 1+t\right) \right) }{e_{q}\left( t\right) -u}e_{\lambda }^{\left(
x_{1}\right) }\left( t\right) e_{\lambda }^{\left( x_{2}\right) }\left(
t\right) \\
&&\left\{ \sin _{\lambda }^{\left( x_{1}\right) }\left( t\right) \cos
_{\lambda }^{\left( x_{2}\right) }\left( t\right) +\cos _{\lambda }^{\left(
y_{1}\right) }\left( t\right) \sin _{\lambda }^{\left( y_{2}\right) }\left(
t\right) \right\}
\end{eqnarray*}%
\begin{eqnarray}
&=&\sum\limits_{m=0}^{\infty }FG_{m,\lambda }^{\left[ k,s\right] }\left(
x_{1},y_{1};u\right) \frac{t^{m}}{m!}\sum\limits_{k=0}^{\infty }C_{k,\lambda
}\left( x_{2},y_{2}\right) \frac{t^{k}}{k!}  \notag \\
&&+\sum\limits_{m=0}^{\infty }FG_{m,\lambda }^{\left[ k,s\right] }\left(
x_{1},y_{1};u\right) \frac{t^{m}}{m!}\sum\limits_{k=0}^{\infty }S_{k,\lambda
}\left( x_{2},y_{2}\right) \frac{t^{k}}{k!}\text{.}  \label{48}
\end{eqnarray}

Using Cauchy product (\ref{47}) and (\ref{48}), we have the following
theorem.

\begin{theorem}
The following relations hold true:%
\begin{eqnarray}
&&FG_{n,\lambda }^{\left[ k,c\right] }\left( x_{1}+x_{2},y_{1}+y_{2};u\right)
\notag \\
&=&\sum\limits_{k=0}^{n}\binom{n}{k}\left\{ FG_{n-k,\lambda }^{\left[ k,c%
\right] }\left( x_{1},y_{1};u\right) C_{k,\lambda }\left( x_{2},y_{2}\right)
-FG_{n-k,\lambda }^{\left[ k,s\right] }\left( x_{1},y_{1};u\right)
S_{k,\lambda }\left( x_{2},y_{2}\right) \right\}  \label{49}
\end{eqnarray}%
and%
\begin{eqnarray}
&&FG_{n,\lambda }^{\left[ k,s\right] }\left( x_{1}+x_{2},y_{1}+y_{2};u\right)
\notag \\
&=&\sum\limits_{k=0}^{n}\binom{n}{k}\left\{ FG_{n-k,\lambda }^{\left[ k,s%
\right] }\left( x_{1},y_{1};u\right) C_{k,\lambda }\left( x_{2},y_{2}\right)
+FG_{n-k,\lambda }^{\left[ k,c\right] }\left( x_{1},y_{1};u\right)
S_{k,\lambda }\left( x_{2},y_{2}\right) \right\} \text{.}  \label{50}
\end{eqnarray}
\end{theorem}

Putting $x_{1}=x_{2}=x$ and $y_{1}=y_{2}=y$ in (\ref{49}) and (\ref{50}),
respectively, we have%
\begin{eqnarray*}
&&FG_{n,\lambda }^{\left[ k,c\right] }\left( 2x,2y;u\right) \\
&=&\sum\limits_{k=0}^{n}\binom{n}{k}\left\{ FG_{n-k,\lambda }^{\left[ k,c%
\right] }\left( x,y;u\right) C_{k,\lambda }\left( x,y\right)
-FG_{n-k,\lambda }^{\left[ k,s\right] }\left( x,y;u\right) S_{k,\lambda
}\left( x,y\right) \right\}
\end{eqnarray*}%
and%
\begin{eqnarray*}
&&FG_{n,\lambda }^{\left[ k,s\right] }\left( 2x,2y;u\right) \\
&=&\sum\limits_{k=0}^{n}\binom{n}{k}\left\{ FG_{n-k,\lambda }^{\left[ k,s%
\right] }\left( x,y;u\right) C_{k,\lambda }\left( x,y\right)
+FG_{n-k,\lambda }^{\left[ k,c\right] }\left( x,y;u\right) S_{k,\lambda
}\left( x,y\right) \right\} \text{.}
\end{eqnarray*}

\section{Conclusion}

Kim and Kim \cite{7} considered the polyexponential and unipoly functions.
Kim et al. (\cite{3} and \cite{11}) defined and investigated the new type
modified degenerate polyexponential function, the degenerate poly-Bernoulli
polynomials and the degenerate poly-Genocchi polynomials. Motivated by these
studying, we introduce the degenerate poly-Frobenius-Genocchi polynomials of
the complex variables. We also give their some interesting properties and
identities. As one of our future projects, we would like to continue to do
researcher on degenerate versions of various special numbers and polynomials.

\begin{acknowledgement}
The present investigation was supported, in part, by the Scientific Research
Project Administration of the University of Akdeniz.
\end{acknowledgement}

\end{document}